\documentclass[11pt]{article}
\usepackage{amscd, amsmath, amssymb}

\title{Curvature, diameter, and quotient manifolds}
\author{Burt Totaro}
\date{  }

\def\Z{\text{\bf Z}}
\def\Q{\text{\bf Q}}
\def\R{\text{\bf R}}
\def\C{\text{\bf C}}
\def\P{\text{\bf P}}
\def\arrow{\rightarrow}

\def\surj{\twoheadrightarrow}
\def\qed{\ QED }

\def\ord{\text{ord}}

\begin{document}
\maketitle

\newtheorem{theorem}{Theorem}[section]
\newtheorem{corollary}[theorem]{Corollary}
\newtheorem{lemma}[theorem]{Lemma}

This paper gives improved counterexamples
to a question
by Grove (\cite{Grove}, 5.7).
The question was whether for each positive integer
$n$ and
real number $D$, the simply connected closed Riemannian $n$-manifolds $M$ with
sectional
curvature $\geq -1$ and diameter $\leq D$ fall into only finitely many
rational homotopy types. This was suggested by Gromov's theorem which
bounds the Betti numbers of $M$ in terms of $n$ and $D$ \cite{Gromov}.
It was known that there can be infinitely many integral homotopy types
already in dimension 7, perhaps first by Aloff and Wallach \cite{AW}.

Fang and Rong recently gave a negative answer to Grove's question
in all dimensions
$\geq 22$ (\cite{FR}, Theorem B). We use certain biquotient manifolds,
that is, quotients of homogeneous manifolds $G/H$ by a subgroup
of $G$ which acts freely,
to show that
the question has a negative answer already in dimension 6. Our examples are
in fact nonnegatively curved. More precisely,
we find infinitely many rational cohomology rings among simply
connected closed Riemannian 6-manifolds with nonnegative
sectional curvature. (Of course, we can arrange that these manifolds
also have diameter at most 1, by scaling.) The dimension 6 here
is optimal, meaning that
Grove's question has a positive answer in dimensions $\leq 5$.
This follows from Gromov's bound on the Betti numbers, since
the Betti numbers of a simply connected manifold of dimension $\leq 5$
determine its rational homotopy type up to finitely many possibilities.
More precisely, the conjecture that simply connected manifolds of
nonnegative curvature are integrally elliptic would imply, by
Paternain and Petean (\cite{PP}, Corollary 3.6),
that simply connected 5-manifolds of nonnegative curvature
fall into only 4 diffeomorphism classes: $S^5$, $S^3\times S^2$,
the nontrivial $S^3$-bundle over $S^2$, and the Wu manifold
$SU(3)/SO(3)$ \cite{Dold}.

Fang and Rong's examples have the merit of also having an upper
bound on curvature. That is, for $n\geq 22$,
Fang and Rong find numbers $C$ and $D$ such that there are
infinitely many rational
cohomology rings among simply connected closed Riemannian $n$-manifolds
with curvature $-1\leq K \leq C$ and diameter $\leq D$. The next main
result of this paper is that such examples
exist already among $7$-manifolds. This is optimal,
since Tuschmann \cite{Tuschmann} and independently Fang and Rong
\cite{FRcircle} have proved
that in dimensions $\leq 6$ there are only finitely many diffeomorphism
classes in the given class of manifolds. Finally, in dimension 9,
we use biquotients to
give a similar counterexample using only nonnegatively curved manifolds.
That is, for some $C$ and $D$, there are
infinitely many rational cohomology rings among
simply connected closed $9$-manifolds with curvature $0\leq K
\leq C$ and diameter $\leq D$.

To conclude, one can ask what substitute for Grove's question
might be true. For the problem with an upper curvature bound,
there is already a remarkable substitute for Grove's question,
the Petrunin-Tuschmann theorem (\cite{PT}, Corollary 0.2).
Namely, for each $n$,
$C$, and $D$, there is a finite set of closed smooth manifolds $E_i$
of dimension $\geq n$
such that any simply connected closed Riemannian $n$-manifold
with curvature $-1\leq K\leq C$ and diameter $\leq D$ is diffeomorphic
to the quotient of some $E_i$ by a free action of a torus. Thus,
in Fang and Rong's examples and in our
examples in dimensions 8 and 9, the infiniteness
comes entirely from considering quotients of a single manifold
by different free torus actions.

In section \ref{conj}, we suggest some possible substitutes for
Grove's question with no upper curvature bound.

I would like to thank Gabriel Paternain, Xiaochun Rong, and Krishnan
Shankar for their papers, as well as many useful conversations.

\section{Counterexamples to Grove's question among nonnegatively
curved 6-manifolds}

Here we prove:

\begin{theorem}
\label{6}
There are infinitely many
isomorphism classes of
rational cohomology rings among
simply connected closed Riemannian $6$-manifolds
with nonnegative sectional curvature.
\end{theorem}

As explained in the introduction, this gives a negative answer to
Grove's question
in dimension 6, which is optimal. Also, it follows from the
theorem of Tuschmann \cite{Tuschmann} and Fang-Rong \cite{FRcircle}
discussed in the introduction
that there cannot be an upper bound on the curvature of the manifolds
we construct, if we fix their diameter to be 1.

{\bf Proof. }The 6-manifolds $M$ we construct will all be biquotients,
of the form $(S^3)^3/(S^1)^3$ for different free isometric actions
of the group $(S^1)^3$ on the Riemannian manifold $(S^3)^3$.
Like all biquotients, these manifolds $M$ have nonnegative sectional
curvature, by O'Neill's curvature formula for Riemannian
submersions \cite{O'Neill}, which we state in the proof of 
Theorem \ref{7}.

Let $(S^1)^2\subset SO(4)$ be the standard maximal torus, with
$(\lambda,\mu)\in (S^1)^2$ acting isometrically on $S^3\subset \C^2$
by $(\lambda,\mu)(u,v)=(\lambda u,\mu v)$. Therefore we  have a natural
isometric action of $(S^1)^6$ on $(S^3)^3$. The actions of $(S^1)^3$
on $(S^3)^3$ we consider will be given by homomorphisms
$(S^1)^3\arrow (S^1)^6$, which we will specify further as we go along.
As a first simplification, let us assume that the homomorphism
$(S^1)^3\arrow (S^1)^6\surj (S^1)^3$ which gives the action of
$(S^1)^3$ on the coordinates $(u_1,u_2,u_3)$
is the identity. Given this, the homomorphisms $(S^1)^3\arrow (S^1)^6$
we consider will be determined by a $3\times 3$ matrix of integers
$$\begin{pmatrix} a_1 & a_2 & a_3 \\ b_1 & b_2 & b_3 \\
c_1 & c_2 & c_3 \end{pmatrix},$$
with the action of $(S^1)^3$ on $(S^3)^3$ given by
\begin{multline*}
(\lambda_1,\lambda_2,\lambda_3)((u_1,v_1),(u_2,v_2),(u_3,v_3)) \\
=((\lambda_1 u_1, \lambda_1^{a_1}\lambda_2^{a_2}\lambda_3^{a_3}v_1),
(\lambda_2 u_2, \lambda_1^{b_1}\lambda_2^{b_2}\lambda_3^{b_3}v_2),
(\lambda_3 u_3, \lambda_1^{c_1}\lambda_2^{c_2}\lambda_3^{c_3}v_3)).
\end{multline*}

To check whether this action of $(S^1)^3$ on $(S^3)^3$ is free, 
one sees easily that it
suffices to check freeness at the 8 points $(p_1,p_2,p_3)$ in $(S^3)^3$
with each $p_i$ equal to $(1,0)$ or $(0,1)$.
Here freeness
at the point $((1,0),(1,0),(1,0))$ is automatic by our choice
of the action on the variables $u_i$. Freeness at the 3 points
with  one $(0,1)$ component means that the diagonal entries
$a_1,b_2,c_3$ of our matrix are $\pm 1$. Freeness at the 3 points
with two $(0,1)$ component means that the 3 determinants
$$\det\begin{pmatrix}b_2 & b_3\\ c_2 & c_3\end{pmatrix},
\det\begin{pmatrix}a_1 & a_3\\ c_1 & c_3\end{pmatrix},
\det\begin{pmatrix}a_1 & a_2\\ b_1 & b_2\end{pmatrix}$$
are $\pm 1$. Finally, freeness at the point $((0,1),(0,1),(0,1))$
means that the whole $3\times 3$ matrix has determinant $\pm 1$.

Let us choose our
$3 \times 3$ integer matrix to have the form
$$\begin{pmatrix} 1  & 0 & 0 \\ b_1 & 1 & 1 \\
c_1 & 2 & 1 \end{pmatrix}.$$
Then the above conditions are satisfied for all integers $b_1$ and $c_1$.
Thus the corresponding actions of $(S^1)^3$ on $(S^3)^3$ are free.
We will show that the corresponding quotient manifolds
$M=(S^3)^3/(S^1)^3$ have infinitely many non-isomorphic rational
cohomology rings.

To compute the cohomology ring of such a  manifold $M$, we consider
the associated fibration
$$(S^3)^3\arrow (S^3)^3/(S^1)^3\arrow (BS^1)^3.$$
Thus we can consider $M$, up to homotopy, as being obtained
from $(BS^1)^3$ by passing three times from the base space to the 
total space of an $S^3$-bundle. The cohomology ring of
$(BS^1)^3=(\C\P^{\infty})^3$ is the polynomial ring
$\Z[x_1,x_2,x_3]$. By our choice of the action of $(S^1)^3$ on
$(S^3)^3$, our three $S^3$-bundles over
$(BS^1)^3$ have Euler classes in $H^4((BS^1)^3,\Z)$ of the form
$x_1^2$, $x_2(b_1x_1+x_2+x_3)$, and $x_3(c_1+2x_2+x_3)$. We observe that
these elements form a regular sequence in the polynomial
ring $\Z[x_1,x_2,x_3]$. Therefore, applying the spectral sequence
successively for these three $S^3$-bundles shows that
$$H^*(M,\Z)=\Z[x_1,x_2,x_3]/(x_1^2, x_2(b_1x_1+x_2+x_3),
x_3(c_1+2x_2+x_3)).$$

It remains to show that these 6-manifolds $M$ have infinitely many
non-isomorphic rational cohomology rings, as the integers $b_1$
and $c_1$ vary. It turns out that later calculations are slightly
simplified if we define the rational numbers
$a=c_1/4$ and $b=(2b_1-c_1)/4$, so that
$b_1=2(a+b)$ and $c_1=4a$. Things will also simplify if we multiply the second
relation by 2. In these terms,
$M$ has rational cohomology ring
$$H^*(M,\Q)=\Q[x_1,x_2,x_3]/(x_1^2,2x_2(2(a+b)x_1+x_2+x_3),
x_3(4ax_1+2x_2+x_3)).$$
We will only consider the manifolds $M$ associated to integers
$b_1$ and $c_1$ which are not both zero. Then $a$ and $b$ are not both zero.

It turns out that we only need to consider the rational cohomology ring
of $M$ in degrees $\leq 4$. In this range, the ring is described
by a 3-dimensional vector space $V=H^2(M,\Q)$ together with a
3-dimensional linear subspace of $S^2V$, the kernel of the product
map $S^2V\arrow H^4(M,\Q)$.
We need to extract a more understandable invariant from this 3-dimensional
linear system of quadrics, which is spanned by the above 3
relations. One approach, which I learned from Wall \cite{Wallnet},
is to consider the determinant of the
quadratic form given by an arbitrary linear combination of
the above 3 relations, 
$$\lambda x_1^2+\mu \cdot 2x_2(2(a+b)x_1+x_2+x_3)
+\nu x_3(4ax_1+2x_2+x_3).$$
The determinant of this quadratic form is:
$$\det\begin{pmatrix}\lambda & 2(a+b)\mu & 2a\nu\\ 2(a+b)\mu & 2\mu & 
\mu+\nu\\ 2a\nu & \mu+\nu & \nu \end{pmatrix}=
-\lambda\mu^2-\lambda\nu^2+4(a^2-b^2)\mu^2\nu+8ab \mu\nu^2.$$
This is a more geometrically understandable object, a cubic curve
over $\Q$ which is an invariant of the rational cohomology ring
of $M$, modulo the action of $GL(3,\Q)$ on $\lambda,\mu,\nu$
and modulo scalars. For convenience, write $\alpha=4(a^2-b^2)$
and $\beta=8ab$, so the cubic has the form
$$-\lambda\mu^2-\lambda\nu^2+\alpha\mu^2\nu+\beta \mu\nu^2=0.$$
Since $a$ and $b$ are not both zero,
$\alpha$ and $\beta$ are not both zero.
Then we compute that the cubic curve has exactly one singular point,
a node at the point $[1,0,0]$ in $\P^2$.

We still need to extract a more computable invariant from this
nodal cubic over $\Q$. We use that this curve has 3 inflection points
over the algebraic closure of $\Q$, apart from the singular point.
The lines from these 3 inflection points to the singular point
$[1,0,0]$ are specified by the point $[\mu,\nu]\in\P^1$ associated
to each inflection point $[\lambda,\mu,\nu]\in\P^2$. Here we are
thinking of $\P^1$ as the space of lines through the singular
point $[1,0,0]$ in $\P^2$. Computing
shows that lines through the 3 inflection points
are the 3 roots of the equation:
$$\beta \mu^3-3\alpha \mu^2\nu-3\beta \mu\nu^2+\alpha \nu^3=0.$$

The two tangent lines to the nodal cubic curve at its singular point
are described by the binary quadratic form $\mu^2+\nu^2$.
Therefore, knowing the nodal cubic curve
modulo scalars and automorphisms of $\P^2$ determines
the above binary cubic form modulo scalars
and modulo automorphisms of $\P^1$
which preserve the binary quadratic form $\mu^2+\nu^2$ up to scalars.
This automorphism group is an orthogonal group $O(2)$ times
the scalars, and so it has two connected components. The identity
component consists of the automorphisms
\begin{align*}
\mu &\mapsto c\mu+d\nu \\
\nu &\mapsto -d\mu+c\nu
\end{align*}
with $c,d\in \Q$, not both zero.
We compute that this group acts on the above binary cubic by the
following change in $\alpha$ and $\beta$:
$$\alpha+\beta i\mapsto (\alpha+\beta i)(c+di)^3,$$
where $i=\sqrt{-1}$ as usual. All our notation has been chosen in order
to make this formula as simple as possible.

Thus, to the above binary cubic we associate the number
$\alpha+\beta i$ in $K=\Q(i)$. It is nonzero since we assumed
$\alpha$ and $\beta$ are not both 0. Scaling the binary cubic leaves
a well-defined class in $K^*/\Q^*$. Making the above coordinate change
leaves a well-defined class in $K^*/(\Q^*\cdot (K^*)^3))=(K^*/3)/(\Q^*/3)$.
Finally, we have to consider the effect of a coordinate change not
in the identity component of the above orthogonal group; for example,
we can switch $\mu$ and $\nu$. This changes $\alpha+\beta i$
to $\beta+\alpha i$. To sum up, we can say that the rational cohomology
ring of the 6-manifold $M$ determines an unordered pair
of elements of the group $(K^*/3)/(\Q^*/3)$, the class of
$\alpha+\beta i$ and the class of $\beta+\alpha i$.

We can now see that we have infinitely many isomorphism
classes of rational cohomology rings among the 6-manifolds $M$.
First, the group $(K^*/3)/(\Q^*/3)$ is infinite. This follows,
for example, from the existence of infinitely many prime numbers $p$
which split in $K=\Q(i)$, namely all primes $p\equiv 1\pmod{4}$.
If $\pi_1$ and $\pi_2$ are the two prime ideals in $K$ which lie over $p$,
then $\ord_{\pi_1}(x)-\ord_{\pi_2}(x)$ is a surjective homomorphism
from $(K^*/3)/(\Q^*/3)$ to $\Z/3$. Since we have infinitely many
such homomorphisms, only finitely many of which can be nontrivial
on a given element of $K^*$, the group $(K^*/3)/(\Q^*/3)$ is infinite.

Furthermore, any element of the group $(K^*/3)/(\Q^*/3)$
has the form $\alpha+\beta i$ for some $\alpha,\beta$ coming
from a 6-manifold $M$ as above. Indeed, the 6-manifold $M$
is described by a pair of integers $b_1,c_1$, and then our definitions
say that
\begin{align*}
a &= c_1/4\\
b &= (2b_1-c_1)/4\\
\alpha+\beta i &= 4(a+bi)^2.
\end{align*}
Since the group $(K^*/3)/(\Q^*/3)$ is 3-torsion, every element $\alpha
+\beta i$ is a square in this group, and so every element can be written as
$4(a+bi)^2$ for some $a,b\in \Q$. Also, multiplying $a+bi$ by
any nonzero integer does not change its class in $(K^*/3)/(\Q^*/3)$,
so every element of the latter group has the form $4(a+bi)^2$ for
some $a,b\in \Z$. Then the corresponding numbers $b_1$ and $c_1$
are also integers. Thus we have
shown that the 6-manifolds we consider can give rise to any element
of the infinite group $(K^*/3)/(\Q^*/3)$.

To be precise, the invariant
of the rational cohomology ring we defined is an unordered pair of
elements of this group. This is enough to show that there are infinitely
many isomorphism clases of rational cohomology rings among these
6-manifolds. \qed

\section{Counterexamples to Grove's question among 7-manifolds with
upper curvature bound}

Here we prove:

\begin{theorem}
\label{7}
There are numbers $C$ and $D$ such that there are infinitely many
isomorphism classes of
rational cohomology rings among
simply connected closed Riemannian $7$-manifolds
with curvature $-1\leq K \leq C$ and diameter $\leq D$.
\end{theorem}

As mentioned in the introduction, this strengthens the examples
of Fang and Rong (\cite{FR}, Theorem B), by lowering the
dimension from 22 to 7. The dimension 7 is optimal, by the
theorem of Tuschmann  \cite{Tuschmann} and Fang-Rong \cite{FRcircle},
as discussed in the introduction.

{\bf Proof. }The manifolds we construct will all be quotients
of a fixed manifold $E$ by different free torus actions. In fact,
this is the only way to get examples as in the theorem,
by the Petrunin-Tuschmann theorem
(\cite{PT}, 0.2), as discussed in the introduction.

In this case, $E$ will be an 11-manifold with a free $(S^1)^5$-action.
Before defining $E$, we will construct a 6-manifold $M$ which
will be the quotient $E/(S^1)^5$. By Sullivan \cite{Sullivan}, Theorem 13.2,
for any graded-commutative $\Q$-algebra which is $\Q$ in degree 0,
0 in degree 1, and which satisfies Poincar\'e duality
of dimension 6, there is a smooth simply connected 6-manifold
$M$ with the given rational cohomology ring. (One could also use
Wall's more precise results on 6-manifolds \cite{Wall6}.)
We will take $M$
to have $b_2=5$ and $b_3=0$. Let $x_0,\ldots,x_4$ be a basis for
$V:=H^2(M,\Q)$. Then we choose $M$ so that the cubic form on 
$H^2(M,\Q)$ is given by
$$\int_M
(a_0x_0+\ldots+a_4x_4)^3=c(a_0^2a_1+a_1^2a_2+a_2^2a_3+a_3^2a_4+a_4^2a_0)$$
for some nonzero constant $c$.
The corresponding cubic 3-fold in $\P^4$ is known as the Klein cubic.
It is not the only cubic form one could use,
but calculations with
it are particularly easy. For example, Adler used the Klein cubic as a tool
to describe the Hessian quintic 3-fold of a general cubic 3-fold
(\cite{Adler}, Appendix IV). Using the above formula for the cubic form,
it is easy to check that
the product map $S^2H^2(M,\Q)\arrow H^4(M,\Q)$ is surjective.

Since $M$ is simply connected, $H^2(M,\Z)$ is torsion-free and hence
isomorphic to $\Z^5$, since $H^2(M,\Q)$ is the 5-dimensional vector
space $V$. Let $E$ be the total space of 
the corresponding $(S^1)^5$-bundle over $M$. Thus $E$ is an 11-manifold
with a free $(S^1)^5$-action. Choose a Riemannian metric on $E$
preserved by the torus action. We can scale the metric so as to have
curvature $\geq -1$. Let $D$ be the diameter of $E$.

We consider the 7-manifolds $Y$ which are quotients of $E$
by subtori $(S^1)^4 \subset (S^1)^5$. These
7-manifolds are all simply connected. With the metric induced
from $E$, they all have diameter
$\leq D$. By O'Neill's formula
\cite{O'Neill},
sectional curvature increases under Riemannian submersions, and
so all these quotient manifolds have curvature $\geq -1$. For clarity,
we recall O'Neill's formula here. Let $\pi:E\arrow M$ be a Riemannian
submersion. Let $X,Y$ be linearly independent vector fields on an open
subset of $E$ which are orthogonal to the fibers of $\pi$ (``horizontal''
vector fields). Let $K$ denote sectional curvature on $E$ and $M$. Then
$$K(\pi_*X,\pi_*Y)=K(X,Y)+3|(\nabla_X Y)_v|^2/|X\wedge Y|^2,$$
where the subscript $v$ denotes the projection to the tangent bundle
of the fibers of $\pi$ (``vertical'' projection).

Because $(S^1)^5$ acts freely on $E$, there is also an upper bound
$C$ for the sectional curvature of all quotients of $E$ by subtori
$(S^1)^4\subset (S^1)^5$. The idea here goes back to
Eschenburg \cite{Eschenburg}, Proposition 22. Namely, O'Neill's formula
shows that the sectional curvature of a quotient manifold $E/(S^1)^4$
can be computed locally on $E$. Furthermore, the same formula
for the curvature
formally makes sense for the non-closed subgroup of $(S^1)^5$ associated
to any real linear subspace $\R^4$ in the Lie algebra $\R^5$ of $(S^1)^5$,
using that $(S^1)^5$ acts freely on $E$.
The ``curvature'' so defined is continuous on the compact manifold
of all subspaces $\R^4\subset \R^5$ and all 2-planes in the tangent bundle
of $E$ which are orthogonal to the associated foliation of $E$.
Therefore, there is a uniform upper bound for this curvature function,
and hence for the curvature of all quotients $E/(S^1)^4$ associated
to subtori $(S^1)^4\subset (S^1)^5$.

It remains to show that the different
7-manifolds $Y$ have infinitely many
non-isomorphic rational cohomology rings. We will only consider the
cohomology ring in degrees $\leq 4$. Here $Y$ is an $S^1$-bundle over
the 6-manifold $M$.
We use the corresponding spectral sequence to compute the rational
cohomology of $Y$. First,
$W:=H^2(Y,\Q)=H^2((BS^1)^4,\Q)$
is the 4-dimensional space $V/(\Q\cdot y)$, where $y$ in
$V=H^2((BS^1)^5,\Q)$ is dual to the 
subtorus $(S^1)^4\subset (S^1)^5$ used to define $Y$. Also,
by the same spectral sequence,
the image of $S^2W$ in $H^4(Y,\Q)$ is $H^4(M,\Q)/(y\cdot V)$.

It seems that ``general'' $S^1$-bundles $Y$ over $M$ will not have
interesting cohomology rings. Fortunately, we can exhibit a special
class of $S^1$-bundles which do have interesting cohomology rings.
Namely, consider $S^1$-bundles $Y$ over $M$ corresponding to
elements $y=a_0x_0-a_1x_1+(a_1^3/a_0^2)x_3$ in $V=H^2(M,\Q)$,
for nonzero rational numbers $a_0$ and $a_1$. (For each nonzero
element $y$ up to scalars in $H^2(M,\Q)$, there is a corresponding
$S^1$-bundle $Y$ over $M$ which is simply connected.) That is, $y$
is a rational point on a certain cuspidal cubic curve in the projective
space $\P^4$ of lines in $V$. 
We will see that the corresponding
$S^1$-bundles $Y$ over $M$ have infinitely many non-isomorphic
rational cohomology rings.

The first useful property of points $y$ as above is that multiplication
by $y$, from $V=H^2(M,\Q)$ to $H^4(M,\Q)\cong V^*$, is not an
isomorphism, as we compute directly from the cubic form. In the
19th-century terminology, which we will not really need, this means that
the above cuspidal curve lies on the Hessian quintic 3-fold of the
given cubic 3-fold. In fact, we compute that the kernel
of multiplication by $y$, from $H^2M$ to $H^4M$, is 1-dimensional,
spanned by
$$z:=a_0x_0+a_1x_1+(a_0^2/a_1)x_2.$$
I have the impression that what makes the following proof work is
the curious fact that the birational involution of the Hessian quintic
which takes $y$ to $z$ transforms the cuspidal cubic curve of points $y$
to the smooth conic curve of points $z$.
In general, it is well known
that a quintic 3-fold will have many rational curves, but the
difference between these two curves still seems surprising.
In any case, the following proof will not use these vague ideas.

Because multiplication by $y$ from $H^2M$ to $H^4M$ has 1-dimensional
kernel, the cokernel $H^4M/(y\cdot H^2M)$ is also 1-dimensional.
We identified this cokernel with the image of the product map
$S^2(H^2Y)\arrow H^4Y$.
Thus the cup product on the 7-manifold $Y$ determines a nonzero 
quadratic form on $H^2(Y,\Q)$, well-defined up to scalars. Explicitly,
using that $yz=0\in H^4M$,
this quadratic form is defined,
up to scalars, by
$$(u,v):=\int_M uvz $$
for $u,v\in H^2M/(\Q\cdot y)\cong H^2Y$.

Since $H^2(Y,\Q)$ has dimension 4, which is even, the determinant
in $\Q/(\Q^*)^2$ of a quadratic form on $H^2(Y,\Q)$ is not changed
upon multiplying the quadratic form by a nonzero scalar.
Thus the determinant
of the above quadratic form in $\Q/(\Q^*)^2$ is an invariant of
the rational cohomology ring of $Y$. 

Let $a_2=a_0^2/a_1$, so that $z=a_0x_0+a_1x_1+a_2x_2$.
We compute that the quadratic form $(u,v):=\int_M uvz$ on the
5-dimensional space $H^2M$ is given, up to a constant factor, by the matrix
$$\begin{pmatrix}a_1 & a_0 & 0 & 0 &0\\ a_0 & a_2 & a_1 & 0 & 0\\
0 & a_1 & 0 & a_2 & 0\\ 0 & 0& a_2 & 0 & 0\\ 0 & 0 & 0 & 0 & a_0
\end{pmatrix}.$$
This matrix has determinant 0, because this quadratic form on
$H^2M$ has at least a 1-dimensional kernel spanned by $y$. To compute
the determinant of the resulting quadratic form on $H^2M/(\Q\cdot y)
\cong H^2Y$,
we can use the lower right $4\times 4$ minor, which has determinant equal to:
$$-a_2^3a_0=-(a_0^2/a_1)^3a_0=-a_0^7/a_1^3\sim -a_0/a_1 \in \Q^*/(\Q^*)^2.$$

Thus, as we vary the nonzero rational numbers $a_0$ and $a_1$,
the determinant of the quadratic form up to scalars on $H^2Y$
can take arbitrary values in the infinite group $\Q^*/(\Q^*)^2$.
It follows that the 7-manifolds $Y$ have infinitely many
non-isomorphic rational cohomology rings. \qed

\section{Counterexamples to Grove's question among nonnegatively curved
9-manifolds with upper curvature bound}

Here we prove:

\begin{theorem}
\label{9}
There are numbers $C$ and $D$ such that there are infinitely many
isomorphism classes of
rational cohomology rings among
simply connected closed Riemannian $9$-manifolds
with curvature $0\leq K \leq C$ and diameter $\leq D$.
\end{theorem}

It seems an interesting challenge to find out whether the dimension
here can be improved. Another comment is that, at least in slightly
higher dimensions, there are examples as in Theorem \ref{9}
which have infinitely many isomorphism classes of cohomology rings
with complex coefficients.

{\bf Proof. }The manifolds we construct will all be quotients
of a fixed manifold by different free torus actions,
as in Theorem \ref{7}. In fact,
this is the only way to get examples as in the theorem,
by the Petrunin-Tuschmann theorem
(\cite{PT}, 0.2), as discussed in the introduction.

Precisely, we consider biquotients of the form $(S^3)^4/(S^1)^3$,
for the different subgroups $(S^1)^3\subset (S^1)^4$, where
we will specify a free isometric action of $(S^1)^4$ on $(S^3)^4$, not the
obvious one. These
9-manifolds are all simply connected. With the metric induced
from the standard metric on $(S^3)^4$, they all have diameter
$\leq D$ where $D$ is the diameter of $(S^3)^4$. By O'Neill's formula
\cite{O'Neill},
sectional curvature increases under Riemannian submersions, and
so all these quotient manifolds have nonnegative curvature. Finally,
by the same argument as in the proof of Theorem \ref{7}, since the whole
group $(S^1)^4$ acts freely on $(S^3)^4$, there is also an upper bound $C$
for the sectional curvature of all quotients of $(S^3)^4$ by subtori
$(S^1)^3\subset (S^1)^4$.

We now explain the free isometric action of $(S^1)^4$ on $(S^3)^4$
which we will use. We think of $S^3$ as the unit sphere in $\C^2$.
For any $4\times 4$ lower-triangular matrix $A=(a_{ij})$ of integers
 with 1's on the diagonal, the following isometric action of $(S^1)^4$
on $(S^3)^4$ is free:
\begin{multline*}
(\lambda_1,\lambda_2,\lambda_3,\lambda_4)((u_1,v_1),(u_2,v_2),(u_3,v_3),
(u_4,v_4)) \\
=((\lambda_1 u_1, (\prod_j \lambda_j^{a_{1j}})v_1),\ldots,
(\lambda_4 u_4, (\prod_j \lambda_j^{a_{4j}})v_4)).
\end{multline*}
Let $M$ be the quotient 8-manifold $(S^3)^4/(S^1)^4$. By viewing
$M$ as the total space of a 4-fold iterated $S^3$-bundle over
$(BS^1)^4$, we find that $M$ has cohomology ring
\begin{multline*}
H^*(M,\Z)=\Z[x_1,x_2,x_3,x_4]/(x_1^2, x_2(a_{21}x_1+x_2),
x_3(a_{31}x_1+a_{32}x_2+x_3), \\
x_4(a_{41}x_1+a_{42}x_2+a_{43}x_3+x_4)).
\end{multline*}

We now specialize the integers $a_{ij}$ to make the cohomology ring
of the 8-manifold $M$ equal to:
$$H^*(M,\Z)=\Z[x_1,x_2,x_3,x_4]/(x_1^2, x_2^2, x_3(x_1+2x_2+x_3),
x_4(x_1+2x_2+x_4)).$$
Changing variables over $\Q$ by $x_3\mapsto x_3-x_1/2-x_2$ and
$x_4\mapsto x_4-x_1/2-x_2$, we find that
$$H^*(M,\Q)=\Q[x_1,x_2,x_3,x_4]/(x_1^2, x_2^2,x_3^2-x_1x_2,x_4^2-x_1x_2).$$

We want to show that the 9-manifolds $Y=(S^3)^4/(S^1)^3$ associated to subtori
$(S^1)^3\subset (S^1)^4$ can have
infinitely many non-isomorphic rational cohomology rings. 
As in our previous examples, we only need to consider the cohomology ring
of $Y$ in degrees $\leq 4$.
Here $Y$ can be
the $S^1$-bundle over $M$ corresponding to any element of $H^2(M,\Z)$ which
generates a summand of $H^2(M,\Z)$. In particular, for each
$a,b,c\in\Q$, there is an $S^1$-bundle $Y$ over $M$ such that
$H^2(Y,\Q)$ is the quotient of $H^2(X,\Q)$ by the line spanned
by $x_4-(ax_1+bx_2+cx_3)$. Then $H^2(Y,\Q)$ is spanned by
$x_1,x_2,x_3$. We compute using the spectral sequence of this
$S^1$-bundle that the kernel of the cup product $S^2H^2(Y,\Q)
\arrow H^4(Y,\Q)$ is the linear system of quadrics spanned by
$$x_1^2,x_2^2,x_3^2-x_1x_2,(ax_1+bx_2+cx_3)^2-x_1x_2.$$
To prove the theorem, it suffices to show that we obtain infinitely
many linear systems of quadrics modulo coordinate changes in
$GL(3,\Q)$, as $a,b,c\in\Q$ vary. Let us assume that $a,b,c$ are all
nonzero, as we are free to.

The idea is to consider what squares of linear forms in $x_1,x_2,x_3$
belong to this
linear system, over a given field $k$ containing $\Q$. It is easy
to check that $x_1^2$ and $x_2^2$ are the only squares of linear forms,
up to scalars, in the span of $x_1^2$, $x_2^2$, and $x_3^2-x_1x_2$.
So suppose that we have a linear form with coefficients in $k$
whose square is a nonzero
multiple of $(ax_1+bx_2+cx_3)^2-x_1x_2$ plus a linear combination
of $x_1^2,x_2^2,x_3^2-x_1x_2$. By considering the coefficients
of $x_1x_3$ and $x_2x_3$, we see that the given linear form must be,
after multiplying by a constant in $k$, of the form $ax_1+bx_2+tx_3$
for some $t$ in $k$.

So we have to work out for which values of $t$ does the square
$(ax_1+bx_2+tx_3)^2$ belong to our linear system of quadrics.
We compute that, modulo this linear system:
\begin{align*}
(ax_1+bx_2+tx_3)^2 &\equiv (ax_1+bx_2+tx_3)^2-(t/c)[(ax_1+bx_2+cx_3)^2
-x_1x_2] \\
&\equiv (t^2-ct)x_3^2+(1/c)(2abc-2abt+t)x_1x_2  \\
&\equiv (t^2+(1/c)(-c^2-2ab+1)t+2ab)x_1x_2.
\end{align*}
It is easy to check that
$x_1x_2$ is not zero modulo our linear system of quadrics.
So the square $(ax_1+bx_2+tx_3)^2$ belongs to our linear system
of quadrics if and only $t$ satisfies the quadratic equation
$$t^2+(1/c)(-c^2-2ab+1)t+2ab =0.$$

Thus, the given linear system of quadrics contains the squares
of two linear forms
$x_1^2$ and $x_2^2$ defined over $\Q$, together with two others
defined over the quadratic extension of $\Q$ corresponding to the above
equation. It follows that this quadratic extension of $\Q$
is an invariant of the rational cohomology ring of the given
9-manifold $Y$. It remains to show that as the rational numbers $a,b,c$
vary, we obtain infinitely many different quadratic extensions of $\Q$,
and therefore infinitely many isomorphism classes of rational cohomology
rings for these 9-manifolds.

The quadratic extension of $\Q$ given by the above quadratic equation
is specified by the class of its discriminant $\Delta=B^2-4AC$
in $\Q^*/(\Q^*)^2$. Here we have
$$\Delta=4\Big[ \Big(\frac{2ab-c^2-1}{2c}\Big)^2-1\Big].$$
We assume that $\Delta$ is nonzero, as is clearly true for most $a,b,c$.
Since $a,b,c$ can be arbitrary nonzero rational numbers, 
it is easy to see
 that the class of $\Delta$ in $\Q^*/(\Q^*)^2$ can be any element
of the form $4(x^2-1)\sim x^2-1$ for $x\in \Q$, $x\neq \pm 1$.
In particular, $\Delta$
can take infinitely many values in $\Q^*/(\Q^*)^2$: for example, for any
odd prime $p$, we can take $x=p+1$, and then $x^2-1$ has nonzero
image under the homomorphism $\ord_p:\Q^*/(\Q^*)^2\arrow \Z/2$.
Therefore, the 9-manifolds we consider have infinitely many
isomorphism classes of rational cohomology rings. \qed

\section{Some possible substitutes for Grove's question}
\label{conj}

The following questions can be viewed as substitutes for Grove's question.
They are only slight extensions of well-known conjectures.

Bott's conjecture that simply connected
manifolds with nonnegative curvature are elliptic (\cite{FHT}, p.~519)
suggests that
there should be strong restrictions on the homotopy type of manifolds with
curvature $\geq -1$ and diameter $\leq D$. On the other hand,
any conjecture must be compatible with Grove and Ziller's
examples of nonnegatively curved manifolds, including
all $S^2$-bundles over $S^4$ and all $S^3$-bundles over $S^4$ (\cite{GZ},
Theorem B). Any conjecture must also cover the almost nonnegatively curved
manifolds found
by Fukaya and Yamaguchi (\cite{FY}, Theorem 0.18) and Schwachh\"ofer
and Tuschmann \cite{ST}. By definition, a manifold $M$ has almost nonnegative
curvature if for every $\epsilon>0$, $M$ has a Riemannian metric with
$K_M\cdot{\text{diam}}(M)^2>-\epsilon$.
For example, any (linear) sphere bundle over a sphere
has almost nonnegative curvature, by Fukaya and Yamaguchi.

{\bf Question 1. }Is every closed simply connected manifold with nonnegative
sectional curvature pure elliptic?

By definition, a manifold $M$ is pure elliptic if there is a minimal model
for the cochain algebra $C^*(M,\Q)$
such that the space of algebra generators $V=V_{\text{ev}}\oplus V_
{\text{odd}}$ is finite-dimensional (this says that $M$ is elliptic),
while $d(V_{\text{ev}})=0$ and $d(V_{\text{odd}})$ is contained in the
subalgebra generated by $V_{\text{ev}}$ (\cite{FHT}, p.~435).
For example, every biquotient
manifold is pure elliptic.

{\bf Question 2. }Is every closed simply connected manifold with
almost nonnegative curvature elliptic?

An almost nonnegatively curved manifold need not be pure
elliptic. For example, an $S^5$-bundle over $S^3\times S^3$ with nonzero
Euler class has almost nonnegative curvature by Fukaya and Yamaguchi,
but it is not pure elliptic. It is helpful to observe that
for a pure elliptic space $M$,
the Lie algebra $\pi_{\text{ev}}(\Omega M)\otimes\Q$
is abelian, while for $M$
elliptic, it is only nilpotent. Thus Questions
1 and 2 can be viewed as higher-dimensional analogues of the known results
that a manifold with nonnegative curvature has almost abelian
fundamental group \cite{CG}, whereas a manifold with almost nonnegative
curvature only has almost nilpotent fundamental group \cite{FY}.

{\bf Question 3. }Given $n$ and $D$, is there a finite set of closed
Riemannian orbifolds $B_i$ such that every simply connected closed
Riemannian $n$-manifold
with sectional curvature $\geq -1$ and diameter $\leq D$
fibers over some $B_i$ with fiber almost nonnegatively curved?

Here ``orbifolds'' are allowed to have stabilizer groups equal to
any compact Lie groups, not just finite groups. Question 3 is strongly
suggested by Yamaguchi's Main Theorem (p.~318) and Conjecture (p.~323)
\cite{Yamaguchi}. Questions 2 and 3 together would imply that
simply connected $n$-manifolds with curvature $\geq -1$ and diameter $\leq D$
have only finitely many sequences of
rational homotopy groups $\pi_*(M)\otimes\Q$, using Friedlander
and Halperin's determination of the possible rational homotopy groups
of elliptic spaces (\cite{FHT}, p.~441).

% Omit these bibliography lines if there's no bibliography.

\small \sc DPMMS, Wilberforce Road,
Cambridge CB3 0WB, England.

b.totaro@dpmms.cam.ac.uk
\end{document}